\documentclass[preprint,12pt,dvipsnames]{elsarticle}
     \makeatletter
\def\ps@pprintTitle{%
  \let\@oddhead\@empty
  \let\@evenhead\@empty
  \let\@oddfoot\@empty
  \let\@evenfoot\@empty
}
\makeatother
    \usepackage{amssymb}
     \usepackage{amsthm}
    \usepackage{amsfonts}
    \usepackage{amsmath,amssymb}
    \usepackage[applemac]{inputenc}
    \usepackage{color}
    \usepackage{color, colortbl}
    \usepackage{amsmath}
    \usepackage{amssymb}
    \usepackage{array}
    \usepackage{url}
    \usepackage{bbm}
    \usepackage{graphicx,floatrow}
    \usepackage{subfigure}
\usepackage{booktabs}
\usepackage[all]{xy}
\usepackage{mathrsfs,bbm}
\usepackage[dvipsnames]{xcolor}
    \usepackage{tikz}
    \usetikzlibrary{shapes.geometric, arrows, positioning}
    \usepackage{geometry}
    \usetikzlibrary{arrows}
    \usepackage{amsfonts}
    \usepackage{mathrsfs}
    \usepackage{amsmath,amssymb}
    \usepackage[applemac]{inputenc}
    \usepackage{color}
    \usepackage{amsmath}
    \usepackage{amssymb}
    \usepackage{graphicx}
    \usepackage{tikz}
    \usepackage{adjustbox} 
    \usepackage{tabularx} 

\usepackage{wrapfig}
\usepackage{multicol}
\usepackage{multirow} 
\usepackage{tabularx} 
\usepackage{float} 

\begin{document}
\newcommand{\A}{\mathbb{A}}
\newcommand{\Q}{\mathbb{Q}}
\newcommand{\N}{\mathbb{N}}
\newcommand{\Z}{\mathbb{Z}}
\newcommand{\R}{\mathbb{R}}
\newcommand{\M}{\mathbb{M}}
\newcommand{\tx}[1]{\quad\mbox{#1}\quad}
\newcommand{\Aut}{\mathrm{Aut}}
\newcommand{\Inn}{\mathrm{Inn}}
\newcommand{\Sym}{\mathrm{Sym}}
\newcommand{\LSec}{\mathrm{LSec}}
\newcommand{\RSec}{\mathrm{RSec}}

\begin{frontmatter}
\title{{\bf{ A Psychometric and Practical Comparison of Standard Moodle-Based and STACK-Based Step-by-Step Tests in University Calculus}}\tnoteref{label1}}\tnotetext[label1]{}
\author[nor]{Semen  Bodnarchuk}
\ead{sem_bodn@ukr.net}

\author[nor]{Kateryna . Moskvychova}
\ead{moskvychovakateryna@gmail.com}

\author[nor]{Igor  Orlovsky}
\ead{orlovskyi@matan.kpi.ua}

\author[nor]{Olha  Pelekhata}
\ead{pelehata@gmail.com}

\author[nor]{Olena Tymoshenko}
\ead{olenaty@math.uio.no}
\cortext[cor1]{Corresponding author.}
\address[nor]{``Department of Mathematical Analysis  and Probability Theory, NTUU Igor Sikorsky Kyiv Polytechnic Institute, Ukraine"}


\begin{abstract}
This paper compares two formats of online assessment in university Calculus: a standard Moodle-based step-by-step test and a STACK-based step-by-step test. Both tests assess integration by parts and divide the solution into consecutive response fields, but they differ in their validation mechanisms. The standard test relies on predefined scoring patterns, whereas the STACK-based test uses symbolic validation rules. The comparison is based on final manually verified scoring matrices from two student cohorts and follows a Classical Test Theory framework, including score distributions, reliability estimates, response-field-level indicators, and correlation-based measures. The results show high overall performance and ceiling effects in both formats. However, the STACK-based test required fewer manual corrections, showed higher internal consistency, and produced a more coherent relationship between solution steps and the total score. The findings suggest that symbolic validation can improve scoring robustness in step-by-step mathematical assessment, although careful task design remains necessary to reduce ceiling effects and improve differentiation among higher-performing students.
\end{abstract}
\end{frontmatter}
\section{Introduction}

Distance and blended learning have become a stable part of higher education, especially where face-to-face instruction is disrupted by pandemics, military conflicts, migration, or other crises. The rapid move to online mathematics support during the COVID-19 pandemic showed that digital environments can maintain access to learning, but it also exposed difficulties related to feedback, interaction, and assessment~\cite{mullen2021}. In Ukraine, war-related disruptions have made flexible and reliable assessment practices even more important, particularly in mathematics, where regular practice and teacher support are essential~\cite{unicef2023}.

Mathematics assessment cannot be reduced to checking whether a final answer is correct. Multiple-choice and short-answer questions are useful for some purposes, but they often say little about how a student selected a method, transformed expressions, or completed intermediate stages of a solution. Students may follow a partly correct method, make compensating errors, or reach a correct result without demonstrating the intended reasoning. For this reason, computer-aided assessment in mathematics has long been discussed as a way to assess problem-solving skills and the structure of mathematical reasoning, not only isolated final responses~\cite{beevers2003,greenhow2015}.

Partial credit is central to this problem. In paper-based mathematical assessment, teachers can usually recognise meaningful intermediate work and award credit for partly correct solutions. In computer-based assessment, however, a mathematically valid or partly valid response may be marked as incorrect if it is not represented in the scoring scheme or if its input form differs from the expected one. Incorporating partial credit into computer-aided mathematics assessment is therefore important for bringing automated scoring closer to mathematical marking practice~\cite{ashton2006}.

Step-by-step assessment formats address this limitation by dividing a mathematical problem into consecutive response fields that correspond to selected stages of the solution. Such tests allow the teacher to assess intermediate reasoning as well as the final answer. Previous work on step-by-step tests in university mathematics has shown that structured solution formats can support both assessment and quality analysis of mathematical tasks~\cite{OmOrTy2018}. In Moodle, this approach can be implemented through embedded answer fields and predefined scoring patterns, although this requires the teacher to anticipate acceptable full- and partial-credit responses in advance.

STACK offers a more flexible approach to computer-aided mathematics assessment by using symbolic validation. It relies on computer algebra to interpret mathematical expressions, recognise equivalent forms, and define validation rules for full- and partial-credit responses~\cite{sangwin2008,sangwin2013,SangwinKocher2016}. Recent department-wide use of STACK in undergraduate mathematics also shows that symbolic assessment can be scaled beyond isolated tasks, although it requires careful authoring, testing of validation logic, feedback design, and technical support~\cite{davies2024}.

The use of such systems raises a quality-assurance question. A test may be mathematically well designed but still differ in scoring robustness, internal consistency, diagnostic capacity, and the extent to which automated scoring reflects mathematical correctness. Standards for educational testing emphasise that score interpretation should be supported by evidence of reliability, validity, and appropriate test use~\cite{aera2014}. Psychometric studies of mathematical assessment instruments likewise show that content relevance alone is not sufficient; reliability, item behaviour, and internal structure require empirical examination~\cite{greenhow2015,gleason2019}.

Classical Test Theory provides a practical framework for examining these properties through score distributions, reliability estimates, item-level difficulty and discrimination, and item-total relationships. Item Response Theory offers a complementary approach for larger-scale modelling of item functioning and student ability~\cite{HambletonJones1993}. In the present study, the empirical analysis is based on a CTT-oriented framework because the analysed tests consist of consecutive response fields within a single solution process and because the available samples are more suitable for descriptive psychometric comparison than for full IRT modelling.

This paper compares two assessment formats used in university Calculus courses: a standard Moodle-based step-by-step test and a STACK-based step-by-step test. Both formats assess integration by parts and follow the same sequence of solution steps. Their main difference lies in the validation mechanism: the standard format relies on predefined scoring patterns, whereas the STACK-based format uses symbolic validation rules.

The aim of the paper is to compare these two formats from both psychometric and practical perspectives. The analysis considers scoring robustness, manual review requirements, descriptive score distributions, reliability estimates, response-field-level indicators, item-rest and item-total correlations, and the internal correlation structure of the tests. The contribution of the study lies in its empirical comparison of standard Moodle-based and STACK-based step-by-step tests using final manually verified scoring matrices and a quality-assurance framework adapted to multi-step mathematical assessment.

The paper is organised as follows. Section~\ref{sec:St_of_Prob} describes the problem of declining mathematical preparedness and the limitations of conventional online assessment in mathematics. Section~\ref{sec:SbS&STACK_AF} introduces the two assessment formats considered in the study: standard Moodle-based step-by-step tests and STACK-based step-by-step tests. Section~\ref{sec:Mthds_of_QA} presents the quality-assurance workflow and psychometric indicators used in the analysis. Section~\ref{sec:CompAn42Tests} applies this framework to compare the two formats empirically. Section~\ref{sec:Disc&Conc} discusses the findings, limitations, and directions for future research.

\section{Statement of the Problem and Motivation for Step-by-Step Online Mathematics Assessment}
\label{sec:St_of_Prob}

Many students enter university mathematics courses with uneven or insufficient preparation, a difficulty that has also been discussed in studies of the transition from school to university mathematics~\cite{d, v}. A December 2024 survey conducted with the AdaMatS consortium among $120$ mathematics lecturers from $35$ Ukrainian universities provides local evidence of lecturers' perceptions of this concern~\cite{adamats}. Only $8.4\%$ of respondents estimated that at least $80\%$ of students achieved a sufficient level of mathematical mastery. By contrast, $33.3\%$ believed that this level was reached by approximately $40\%$ of students, and $29.2\%$ estimated that it was reached by only approximately $20\%$. Moreover, $34.2\%$ reported a decline in students' mathematical knowledge, while $48.3\%$ expected further deterioration.

These responses indicate a perceived decline in students' mathematical preparedness and highlight the need for more informative assessment practices. This local evidence is consistent with broader concerns about mathematical achievement reported in international studies and policy reports, including PISA and TIMSS, where mathematical literacy is linked to further study, STEM preparation, and educational quality~\cite{nuno2024,oecd2023,cosgrove2014}.

Conventional online tests provide limited information about mathematical reasoning. A student may select an appropriate method but make an algebraic error, complete several steps correctly but fail at the final stage, or obtain the correct answer without demonstrating the intended procedure~\cite{beevers2003,greenhow2015}. Multiple-choice questions may allow guessing, while short-answer questions may reject mathematically equivalent expressions that do not match predefined answer forms. Partially correct reasoning may therefore remain unrecognised.

This issue may be particularly relevant in distance and blended learning environments, where instructors have fewer opportunities to observe students' work directly. Assessment formats that capture intermediate solution steps can therefore provide useful diagnostic information in addition to overall performance.

Step-by-step assessment divides a mathematical problem into consecutive response fields, allowing intermediate solution steps to be evaluated. The next section compares standard Moodle-based and STACK-based implementations.

\section{Step-by-Step and STACK-Based Assessment Formats}
\label{sec:SbS&STACK_AF}

This study considers two step-by-step assessment formats implemented in Moodle: a standard Moodle-based format and a STACK-based format. Both divide a mathematical solution into consecutive response fields, but they differ in their validation mechanisms. The standard format uses predefined scoring patterns for full- and partial-credit responses, whereas STACK applies symbolic validation rules to recognise mathematically equivalent expressions and assign credit~\cite{sangwin2013,SangwinKocher2016}.

\begin{figure}[H]
  \center{\includegraphics[scale=0.59]{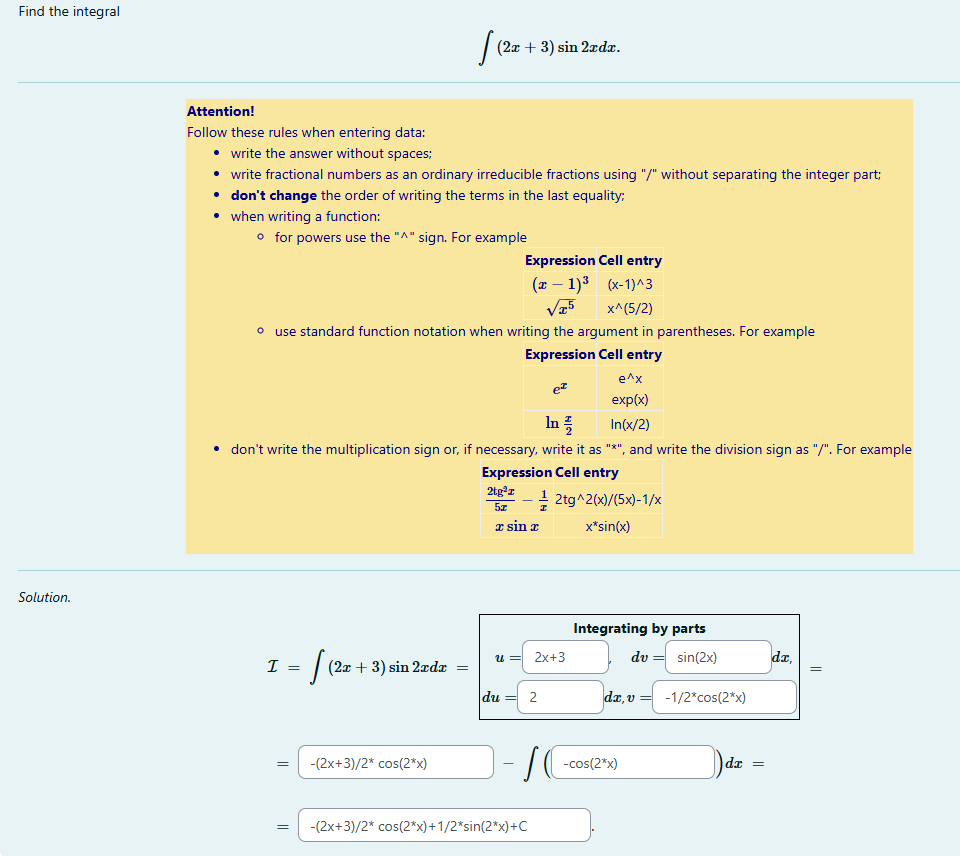}}
 \caption{Example of a standard Moodle-based step-by-step test in which the solution is divided into consecutive completed response fields.}
\label{fig:StandardSbS_test}
\end{figure}

At Igor Sikorsky Kyiv Polytechnic Institute, step-by-step tests have been implemented on the Sikorsky Moodle-based platform for Calculus courses. The examples discussed in this section illustrate the two assessment formats that are later analysed using the quality-assurance framework described in Section~\ref{sec:Mthds_of_QA} and compared empirically in Section~\ref{sec:CompAn42Tests}. The first format is the standard Moodle-based step-by-step test, which uses embedded answer fields and HTML markup to divide the solution into consecutive stages (Figure~\ref{fig:StandardSbS_test}).

\subsection{Standard Moodle-Based Step-by-Step Tests}
\label{subsec:StSbS_Tests}

The standard Moodle-based step-by-step test is designed to make the solution process visible within an online assessment environment. Instead of requiring only a final answer, the task is divided into consecutive response fields corresponding to selected stages of the mathematical solution. In Moodle, this format can be implemented using embedded answer questions and HTML markup, which allow the teacher to combine explanatory text, mathematical notation, and input fields within a single structured task~\cite{OmOrTy2018}.

In the Calculus courses at Igor Sikorsky Kyiv Polytechnic Institute, such tests have been used to assess topics in which the solution can be naturally decomposed into a sequence of mathematical steps. A representative example is integration by parts, where students may be asked to identify $u$ and $dv$, compute $du$ and $v$, construct $uv$ and $vdu$, and evaluate the remaining integral. Figure~\ref{fig:StandardSbS_test} illustrates the standard Moodle-based format in which the solution is divided into consecutive response fields.

In this format, scoring is based on predefined scoring patterns. For each response field, the teacher must specify the forms of responses that receive full or partial credit. This gives the instructor explicit control over the scoring scheme and makes it possible to assess intermediate stages of the solution rather than only the final result. At the same time, the quality of automated scoring depends on how completely mathematically acceptable responses are represented in the predefined patterns.

A practical limitation of this approach is its sensitivity to the form of student input. Mathematically equivalent expressions may fail to receive the intended score if their notation, syntax, or structure does not match one of the predefined scoring patterns. This is particularly relevant for expressions involving products, brackets, fractions, powers, trigonometric functions, or integration notation. As a result, standard step-by-step tests often require clear input instructions for students and periodic review of submitted responses to identify additional full- or partial-credit forms that should be included in the scoring scheme.

\begin{figure}[H]
    \centering
   \includegraphics[width=0.99\textwidth]{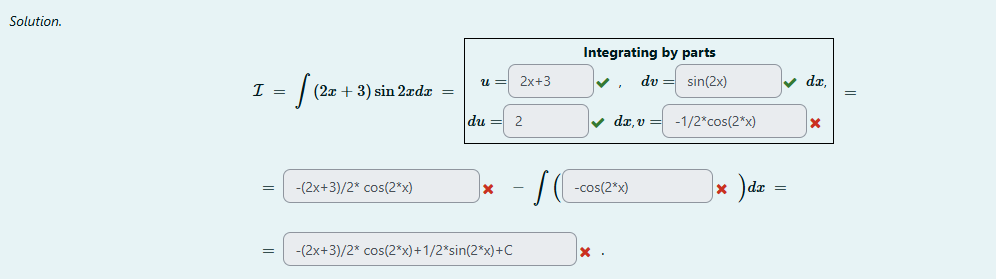}
    \caption{Example of an otherwise mathematically correct response that required manual review due to input formatting}
    \label{fig:partial_error}
\end{figure}

Figure~\ref{fig:partial_error} shows an example of a mathematically acceptable response that requires manual review because its input form differs from the predefined scoring pattern.

\subsection{STACK-Based Step-by-Step Tests and Symbolic Validation}
\label{subsec:STACK_SbST&SV}

The limitations of standard Moodle-based step-by-step tests are mainly connected with their dependence on predefined scoring patterns. STACK offers an alternative approach by integrating computer algebra functionality into Moodle-based assessment. In this format, mathematical expressions entered by students can be interpreted symbolically rather than matched only against predefined textual forms~\cite{sangwin2013,SangwinKocher2016}.

In STACK-based tests, answer validation is organised through symbolic answer tests and validation rules. This makes it possible to recognise mathematically equivalent expressions and to assign full or partial credit according to the mathematical meaning of a response. For step-by-step tasks, this is especially important because each response field may represent a different stage of the solution process. In contrast to standard input matching, the scoring logic does not require the teacher to enumerate every acceptable textual form of an answer.

\begin{figure}[H]    \center{\includegraphics[scale=0.92]{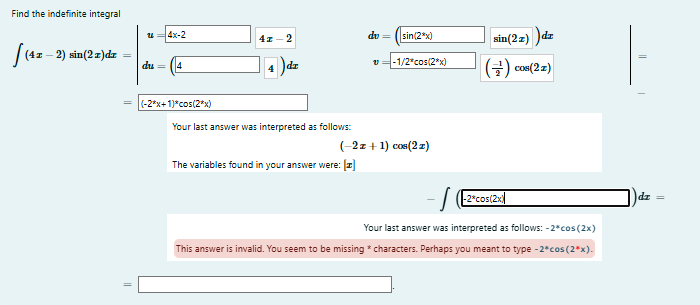}}  
  \caption{Example of a STACK-based step-by-step test with completed response fields and symbolic validation of mathematical expressions.}
  \label{fig:STACK_test}  
\end{figure}

STACK also supports parameterised task templates, which makes it possible to generate comparable variants of the same mathematical problem while preserving the same solution structure. The design of such variants is itself an important aspect of computer-aided mathematics assessment, because variation should preserve the intended mathematical structure rather than introduce unintended differences in difficulty~\cite{SeatonTacy2021}. In the context of integration by parts, the same response-field sequence can be retained: selecting $u$ and $dv$, computing $du$ and $v$, constructing the product term, forming the remaining integral, and entering the final antiderivative. Figure~\ref{fig:STACK_test} illustrates a STACK-based version of such a step-by-step task.

Although symbolic validation reduces dependence on predefined answer patterns, it does not eliminate the need for careful test design. STACK questions require the teacher to define the mathematical structure of the task, specify validation rules, test possible solution paths, and provide clear input guidance for students. Practical reports on STACK development similarly emphasise iterative authoring, testing of validation logic, and refinement of feedback~\cite{zerva2019}. Thus, the advantage of STACK lies not in removing the need for quality assurance, but in providing a more flexible validation mechanism whose empirical performance can be evaluated using the framework described in Section~\ref{sec:Mthds_of_QA}.

\subsection{Comparison of Validation Mechanisms and Implementation Constraints} \label{subsec:Comp_of_VM&IC}

\begin{figure}
    \centering
        \includegraphics[width=0.76\linewidth]{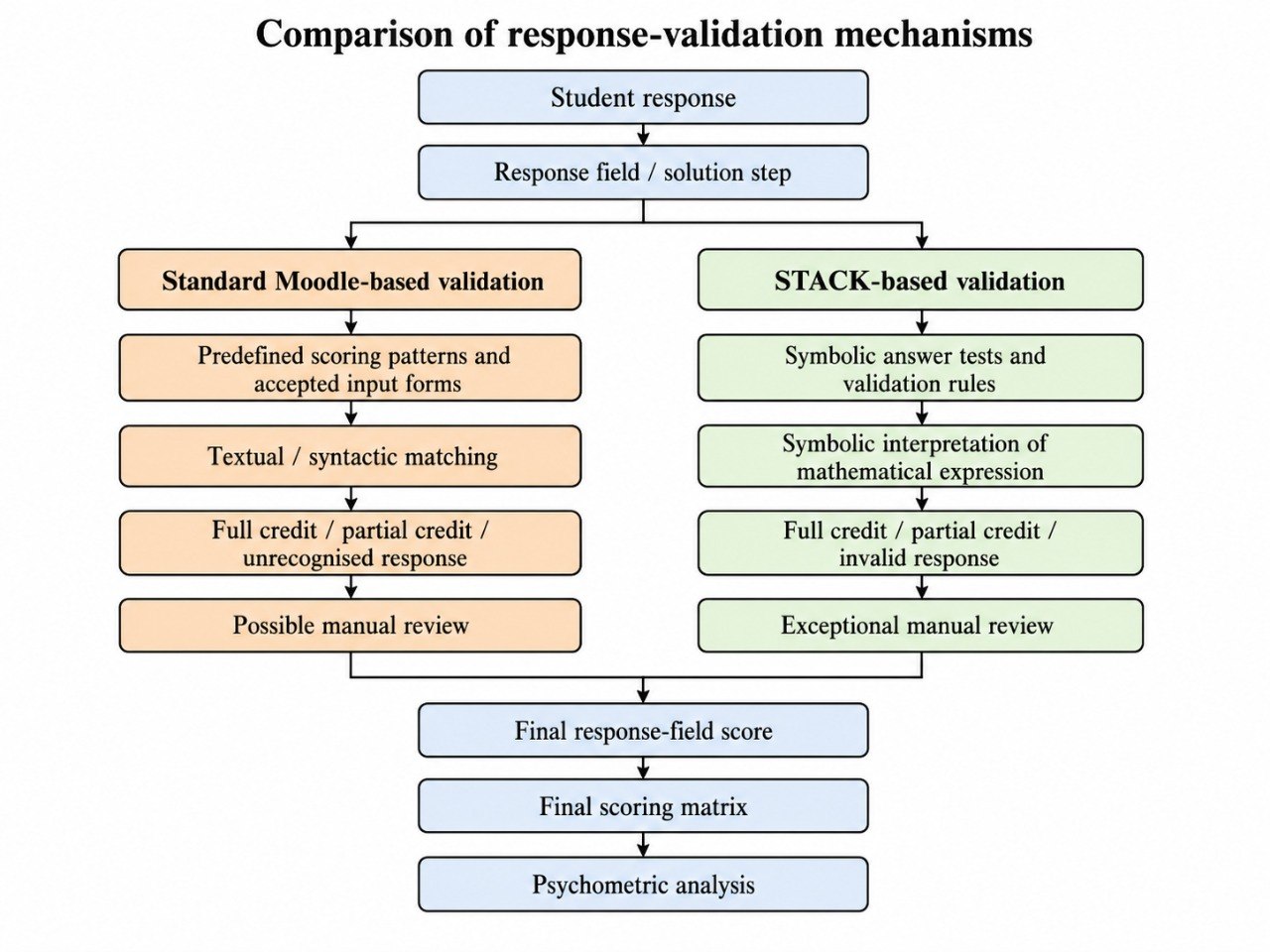}
    \caption{Comparison of response-validation mechanisms in standard Moodle-based and STACK-
based step-by-step tests.}
    \label{fig:valid_compar} 
\end{figure}

The main technical distinction between the two formats lies in how student responses are validated. In the standard Moodle-based format, acceptable responses must be anticipated and encoded as predefined scoring patterns. This approach is relatively transparent for the teacher but requires the manual inclusion of alternative full- and partial-credit forms. In the STACK-based format, validation is based on symbolic answer tests and explicitly defined validation rules, allowing mathematically equivalent expressions to be recognised without enumerating all possible textual forms. Figure~\ref{fig:valid_compar} summarises the technical difference between predefined scoring-pattern validation and symbolic validation.

This difference also affects test maintenance. Standard step-by-step tests may require cumulative expansion of scoring patterns as new mathematically acceptable or partially correct student responses are observed. STACK-based tests require more careful initial authoring, including the definition of variables, validation rules, potential partial-credit cases, and feedback. Thus, the standard format may be easier to implement initially, whereas STACK shifts part of the workload from post-test correction and maintenance to pre-test design and validation.

Both formats therefore require quality assurance, but for different reasons. In standard tests, quality assurance is needed to check whether predefined scoring patterns adequately represent acceptable responses. In STACK-based tests, it is needed to verify whether symbolic validation rules correctly reflect the intended mathematical scoring scheme. These implementation differences motivate the quality-assurance framework described in Section~\ref{sec:Mthds_of_QA} and the empirical comparison presented in Section~\ref{sec:CompAn42Tests}.

\section{ Methods of Quality Assurance} \label{sec:Mthds_of_QA}

\subsection{Quality-assurance workflow and analytical framework}

Quality assurance is essential for determining whether step-by-step mathematics tests provide reliable, interpretable, and diagnostically useful information about students' performance, in line with established standards for educational and psychological testing~\cite{aera2014,CohenSwerdlik2017}. In this study, test quality is considered through three complementary dimensions: reliability, validity understood as alignment between the task content, scoring criteria, and the intended mathematical construct, and diagnostic capacity.

Reliability refers to the consistency with which the set of response fields supports the interpretation of the total test score. In the context of the present study, validity is considered primarily in terms of content alignment: whether the test tasks, solution steps, and scoring rules support the intended interpretation of scores as indicators of students' performance in integration by parts~\cite{aera2014}. Diagnostic capacity refers to the ability of the test to identify differences in students' performance across specific stages of the solution process.

\begin{figure}
    \centering
    \includegraphics[width=0.9\linewidth]{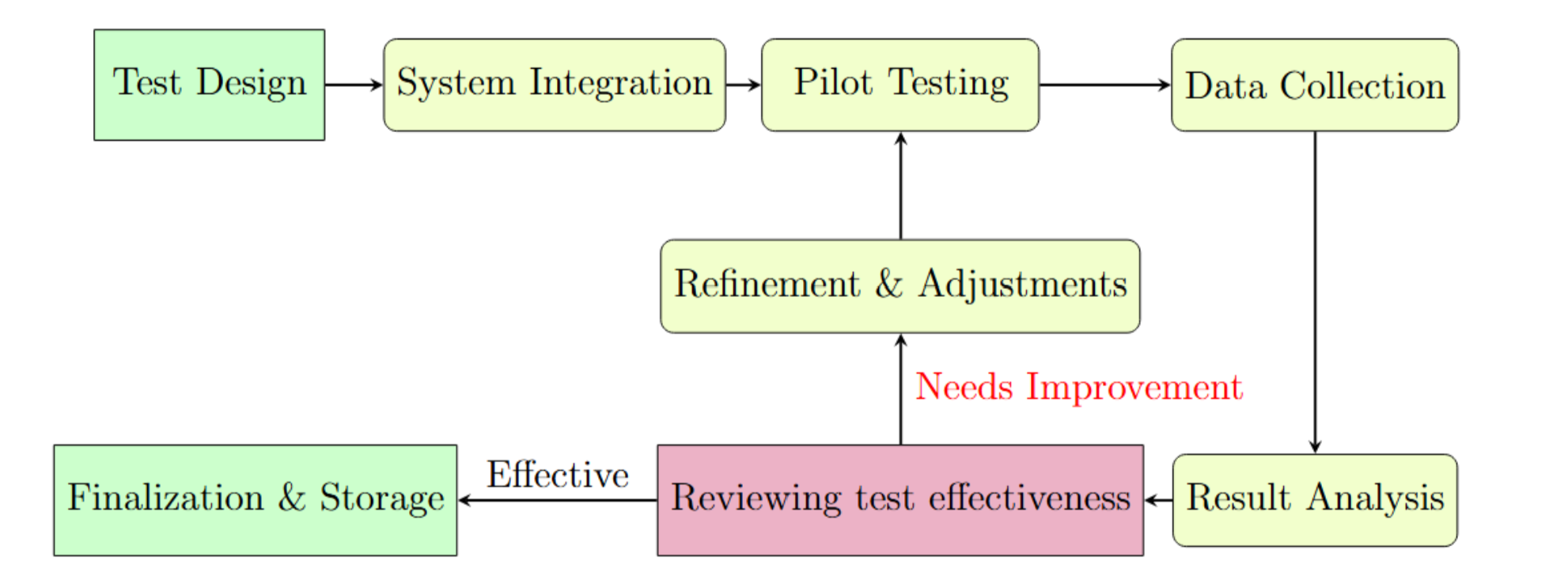}
    \caption{Structured workflow for test quality assurance, outlining test development, data collection, scoring review, statistical analysis, and refinement before further use}
    \label{fig29}
\end{figure}

The quality-assurance workflow used for the analysed tests is shown in Figure~\ref{fig29}. It reflects the quality-control component of computer-based testing guidelines by emphasising test development, pilot testing, data collection, scoring review, statistical analysis, and refinement~\cite{itc2006}. Broader ITC recommendations concerning test administration, security, privacy, and control of delivery are relevant to computer-based assessment but are outside the scope of the psychometric analysis presented here. This workflow is particularly important for step-by-step tests, where scoring depends not only on the final answer but also on intermediate solution steps and on the representation of full- and partial-credit responses.

The empirical analysis follows three main stages. First, student responses are collected from Moodle and organised into a scoring matrix. When automatic scoring produces ambiguous or unrecognised responses, these cases are reviewed manually before the final scoring matrix is constructed. Second, the resulting data are analysed using descriptive statistics, score distributions, reliability coefficients, discrimination indices, and correlation-based item indicators. Third, the results are used to identify solution steps with limited diagnostic value, scoring inconsistencies, ceiling effects, and aspects of the test that require refinement before further use.

Although both Classical Test Theory (CTT) and Item Response Theory (IRT) are commonly used in educational measurement, the present study adopts a CTT-based framework~\cite{HambletonJones1993,fan1998}. This choice is appropriate for the descriptive comparative aim of the study and for the available sample sizes. Because the response fields represent consecutive stages of a single solution process rather than fully independent items, the resulting item-level and correlation-based indicators are interpreted with this structure in mind. IRT-based modelling may provide additional insights in future studies involving larger and more heterogeneous samples. The statistical analysis was carried out using jamovi, Python, and R. 

\subsection{Psychometric indicators used in the analysis}
\label{subsec:PsyInd4An}

The indicators used in this study follow standard approaches to classical test analysis, combining descriptive statistics, response-field-level diagnostics, reliability estimates, and correlation-based indicators~\cite{CohenSwerdlik2017,gulliksen1950,HambletonJones1993}. Similar psychometric approaches have been applied to mathematics tests to examine item difficulty, discrimination, and reliability~\cite{butakor2022}. The psychometric analysis was based on the final scoring matrices constructed after the review of ambiguous or unrecognised responses. Manual review was used to ensure that responses eligible for full or partial credit under the scoring scheme were reflected correctly in the final matrix. Each matrix contained scores for the response fields corresponding to consecutive steps of the solution process. These response fields were treated as scored solution steps rather than as fully independent test items. Scores were normalised to the interval $[0,1]$, and the total test score was computed from the set of response-field scores. This structure made it possible to analyse both overall test performance and the contribution of individual solution steps.

Overall test performance was described using the mean, median, mode, and standard deviation of total scores: the mean characterised the general level of performance, the median and mode indicated score concentration, and the standard deviation reflected score variability. The shape of the score distribution was examined using histograms, boxplots, skewness, and kurtosis. In particular, negative skewness and concentration of scores near the upper end of the scale were interpreted as evidence of a ceiling effect, indicating reduced differentiation among higher-performing students.

At the level of individual response fields, the mean score was used as an indicator of relative difficulty. Higher mean scores indicated that a solution step was completed successfully by most students, whereas lower mean scores reflected greater difficulty or a higher frequency of partial-credit responses. The standard deviation of each response field was used to assess variability in student performance at that stage of the solution.

The discrimination index was used to evaluate how well each response field differentiated between higher- and lower-performing students. It was computed as
\[
D = p_{\mathrm{high}} - p_{\mathrm{low}},
\]
where $p_{\mathrm{high}}$ and $p_{\mathrm{low}}$ denote the mean proportions of the maximum score obtained by the upper- and lower-performing groups, respectively. The use of upper and lower performance groups follows the classical extreme-groups approach to item analysis, applied here to scored response fields~\cite{kelley1939}. As a practical interpretive guideline, values of $D$ around $0.3$ or higher were treated as indicating useful discrimination, whereas lower values suggested limited differentiation between performance levels~\cite{ebel1972}.

Two correlation-based indicators were used to examine the alignment between individual solution steps and the overall score structure. Let $X_j$ denote the score for response field $j$, and let $T=\sum_{k=1}^{m}X_k$ be the total test score across $m$ response fields. The item-total correlation was computed as
\[
\mathrm{ITC}_j = r(X_j, T),
\]
whereas the item-rest correlation was computed as
\[
\mathrm{IRC}_j = r(X_j, T-X_j),
\]
where \(r(\cdot,\cdot)\) denotes the Pearson correlation coefficient. Thus, ITC measures the relationship between a response field and the total score including that field, while IRC relates the response field to the total score computed from the remaining fields and therefore provides a stricter indicator of alignment~\cite{henrysson1963}. In this study, IRC and ITC were interpreted jointly because step-by-step tests consist of consecutive response fields rather than independent items. A response field may discriminate between higher- and lower-performing students while still showing weak alignment with the remaining solution process, especially when later steps depend on earlier responses or when partial-credit scoring is involved.

Internal consistency of the total test score was assessed using Cronbach's $\alpha$ and McDonald's $\omega$. Cronbach's $\alpha$ was computed as
 \[
   \alpha = \frac{m}{m-1} \left(1-\frac{\sum_{j=1}^{m}\sigma_j^2}{\sigma_T^2}\right),
 \]
where $m$ is the number of response fields, $\sigma_j^2$ is the variance of scores for response field $j$, and $\sigma_T^2$ is the variance of the total test score~\cite{cronbach1951}. McDonald's $\omega$ was used as a complementary reliability estimate because it is less dependent on the assumption that all response fields contribute equally to the total score~\cite{mcdonald1999,hayes2020}. It was computed as
 \[
   \omega = \frac{\left(\sum_{j=1}^{m}\lambda_j\right)^2}{\left(\sum_{j=1}^{m}\lambda_j\right)^2+\sum_{j=1}^{m}\psi_j},
 \]
where $\lambda_j$ denotes the factor loading of response field $j$, and $\psi_j$ denotes its error variance. In this study, reliability coefficients were interpreted as indicators of the internal consistency of the total score rather than as evidence that the response fields are independent. Following common practice in educational and psychological measurement, values around $0.7$ or higher were treated as a practical indication of acceptable internal consistency, while lower values suggested that the total score should be interpreted with caution~\cite{CohenSwerdlik2017,nunnally1994}.

Finally, correlation matrices of response fields were constructed to examine the internal structure of each step-by-step test. Let $X_i$ and $X_j$ denote the scores for two response fields. The correlation matrix was defined as
 \[
   R=\Bigl(r(X_i,X_j)\Bigr)_{i,j=1}^m,
 \]
where \(r(\cdot,\cdot)\) denotes the Pearson correlation coefficient computed across students. Because the response fields represent consecutive stages of a single solution rather than independent tasks, the correlations were interpreted as indicators of structural alignment between solution steps. Weak or near-zero correlations were not treated automatically as evidence that a response field was ineffective; they were interpreted in relation to the sequential nature of the task, low score variability, ceiling effects, and possible propagation of errors across later stages.

\section{Comparative Psychometric Analysis of Standard and STACK-Based Step-by-Step Tests}\label{sec:CompAn42Tests}

This section presents a comparative empirical analysis of two assessment formats: a standard step-by-step test and a STACK-based step-by-step test. The comparison focuses on psychometric properties and practical implementation, including scoring robustness, reliability, item-level characteristics, and the implications of the validation mechanism for automated scoring and test maintenance.

The study was conducted at the Faculty of Informatics and Computer Science (FICT) during the autumn semester of the 2024-2025 academic year. Two student cohorts participated: 82 students completed the standard step-by-step test, and 89 students completed the STACK-based test.

As described in Section~\ref{sec:SbS&STACK_AF}, both assessments followed the same step-by-step\linebreak integration-by-parts structure but differed in their validation mechanisms. The analysis applies the indicators defined in Section~\ref{subsec:PsyInd4An} to the final scoring matrices for the two tests. Practical implementation is assessed through the complexity of answer specification, the number of manual score corrections, and the occurrence of input-related evaluation errors.

Because the groups were not randomly assigned, differences in prior mathematical preparation cannot be ruled out. The findings should therefore be interpreted as a comparative description of the two formats rather than as causal evidence of the effects of test format.

The following subsections first examine each test format separately and then compare their psychometric and practical characteristics. Unless stated otherwise, all psychometric analyses in this section were conducted using the final scoring matrices after manual review of responses that were eligible for full or partial credit but were not recognised automatically.

\subsection{Analysis of the Standard Step-by-Step Test}

The standard test used the seven response fields described in Section~\ref{subsec:StSbS_Tests}. Table~\ref{Rslts_SbS} presents an excerpt from the final manually verified scoring matrix. Each row corresponds to a student, and each column contains the score for one solution step. Item scores were normalised to the interval $[0,1]$, and the total test score was calculated as the unweighted mean of the seven item scores.

\begin{table}[H]
    \centering
     \includegraphics[scale=0.8]{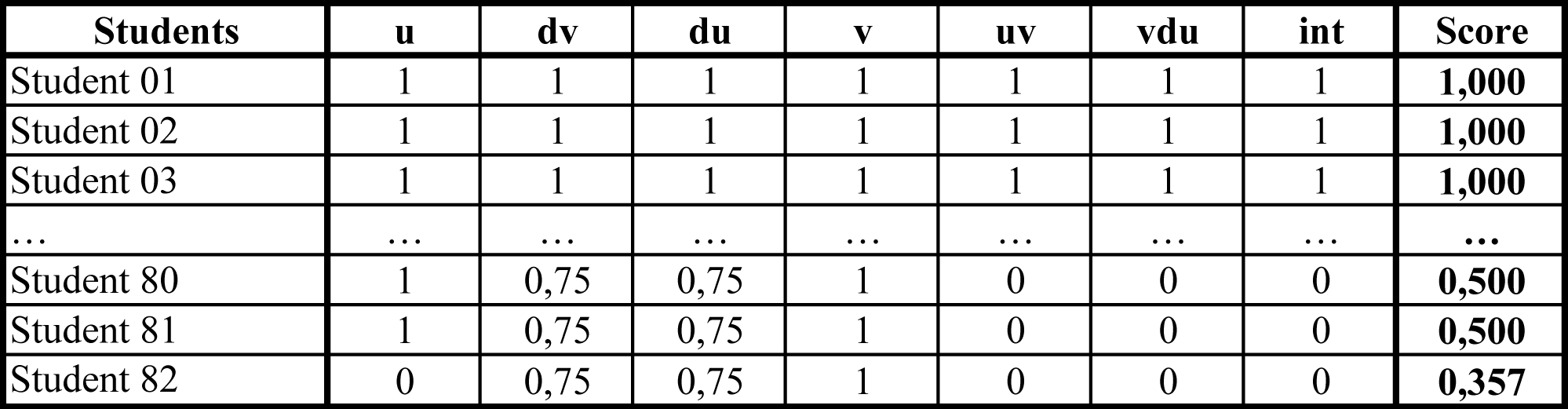}
    \caption{Excerpt from the final response matrix for the standard step-by-step}
   \label{Rslts_SbS}
\end{table}

For the analysed administration, sixteen task variants were prepared manually using a common integration-by-parts template. In each variant, the initial integrand had the form
 \begin{equation}\label{undr_int_fn}
    (Ax+B)\cdot\sin(Cx)\ \ \ \text{or}\ \ \ (Ax+B)\cdot\cos(Cx),
 \end{equation}
where $A$, $B$, and $C$ were integers specified separately for each variant. Thus, all variants required the same mathematical method and sequence of solution steps, while differing in their numerical parameters and trigonometric function.

The standard step-by-step tests analysed in this study have been used and refined over several years. 
For the analysed administration, Table~\ref{tb_SbS_NoV} summarises the number of predefined scoring patterns available for each response field before and after cumulative refinement.

\begin{table}[h]
\centering
\begin{tabular}{l|ccccccc}
\toprule
Characteristics  &  u & dv & du & v & uv & vdu & int \\
\midrule
Initial scoring patterns    & 1 & 2 & 1 & 2 & 3 & 2 & 4  \\
After cumulative refinement  & 2 & 3 & 1 & 5 & 9 & 6 & 18  \\
\bottomrule
\end{tabular}
\caption{Number of predefined scoring patterns for each solution step in the standard step-by-step test}
\label{tb_SbS_NoV}
\end{table}

Table~2 shows that cumulative refinement substantially increased the number of scoring patterns for several response fields, especially \textbf{uv}, \textbf{vdu}, and the final integration step. These increases reflect the diversity of mathematically equivalent full-credit representations and common partially correct forms of student input.

Despite this cumulative refinement, the present administration still required $25$ manual corrections of item scores because some responses that should have received full or partial credit under the scoring scheme were not recognised automatically. As illustrated in Section~\ref{subsec:StSbS_Tests}, this type of mismatch may occur when a mathematically acceptable response differs from the predefined input form. Four additional responses could not be evaluated automatically because of minor syntax deviations; these cases were also reviewed manually and scored according to their mathematical correctness and the corresponding partial-credit rules. The final scoring matrix used in the subsequent analyses incorporates these manual review decisions.

In the present administration, the number of required corrections was treated as an indicator of the practical robustness of the automatic validation mechanism. The psychometric results reported below therefore describe the final, manually verified scoring matrix.

The next part examines the distribution of total scores and the internal-consistency estimates of the standard step-by-step test.

\begin{table}[h]
    \centering
    \begin{tabular}{cccccc}
        \toprule
         \textbf{Mean} & \textbf{Median} & \textbf{Mode} &\textbf{SD} & \textbf{Skewness} & \textbf{Kurtosis} \\
        \midrule
        0.851 & 0.857 & 0.964 & 0.142 & -1.25 & 1.28 \\
        \bottomrule
    \end{tabular}

    \vspace*{1mm}
    \begin{tabular}{cc}
        \toprule
         \textbf{Cronbach's $\alpha$} & \textbf{McDonald's $\omega$}\\
        \midrule
         0.49 & 0.59 \\
        \bottomrule
    \end{tabular}
    \caption{Descriptive statistics and reliability coefficients for the standard test}
    \label{tb_stat_char_SbS}
\end{table}

\begin{figure}[H]
  \center{\includegraphics[scale=0.25]{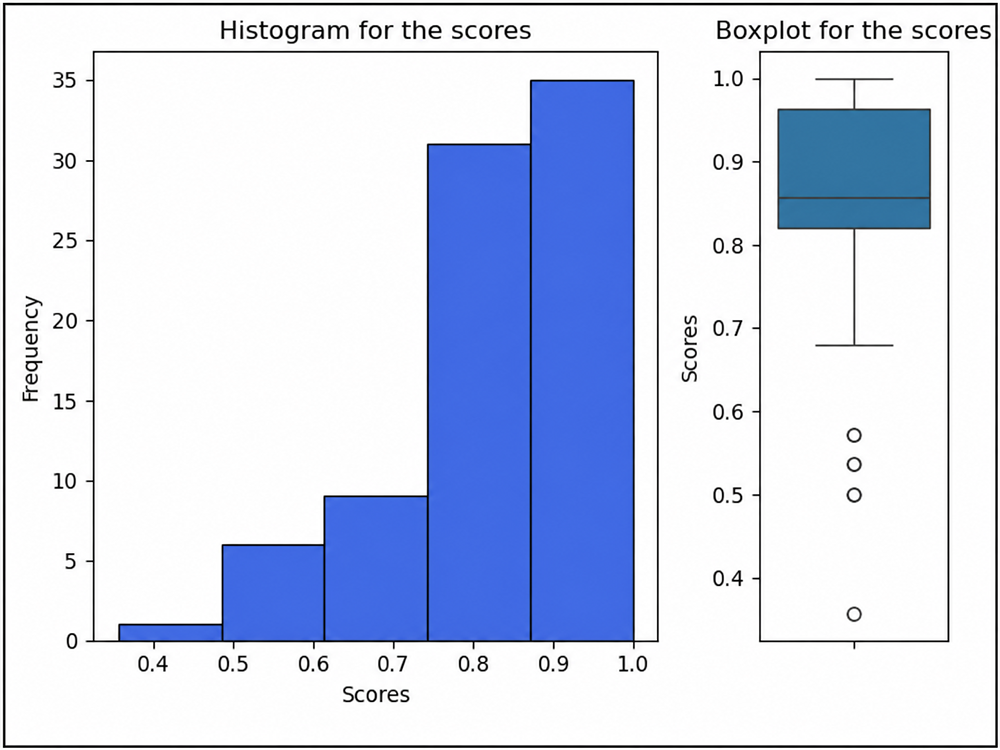}}
  \caption{Histogram and boxplot of the results for the standard step-by-step test}
  \label{fig:SbS_distribution}
 \end{figure}

Based on the final manually verified scoring matrix, the distribution of total scores indicates generally high student performance in the standard step-by-step test. As shown in Table~\ref{tb_stat_char_SbS}, the mean score was $0.851$ and the median was $0.857$, whereas the mode reached $0.964$. Figure~\ref{fig:SbS_distribution} confirms this pattern: the histogram shows a concentration of scores in the upper part of the scale, while the boxplot reveals several lower-end observations. These observations were retained in the analysis, as they reflect genuine variability in student performance.

The concentration of scores near the upper end of the scale suggests a noticeable ceiling effect. Although the test captured some variation in performance, the reduced spread of high scores limited its ability to differentiate among higher-performing students.

The reliability estimates are consistent with this interpretation. Cronbach's $\alpha$ was relatively low ($\alpha = 0.49$), indicating limited internal consistency under the assumption of equal item contributions. McDonald's $\omega$ was higher ($\omega = 0.59$), suggesting a somewhat more favourable but still limited estimate of internal consistency. The difference between the two coefficients is consistent with uneven item behaviour and reduced score variability; in particular, the concentration of high scores weakens inter-item covariance, which directly affects Cronbach's $\alpha$.

These results indicate that the standard step-by-step test adequately captures basic procedural competence but has limited capacity to distinguish among students at the upper end of the performance distribution. This limitation should be considered when interpreting the item-level statistics reported below.

Item-level characteristics of the standard step-by-step test were examined using mean scores, standard deviations, discrimination indices, item-rest correlations (IRC), and item-total correlations (ITC), as reported in Table~\ref{tb_IRC_ITC_SbS}.

\begin{table}[h]
    \centering
    \begin{tabular}{c|ccccccc}
        \toprule
        Response field  &   u  &  dv  &  du  &  v   &  uv  &  vdu & int \\
        \midrule
        Mean Score     & 0.96 & 0.93 & 0.91 & 0.97 & 0.89 & 0.54 & 0.76 \\
        SD             & 0.19 & 0.2  & 0.23 & 0.16 & 0.3  & 0.46 & 0.35 \\
        Discrimination & 0.09 & 0.21 & 0.26 & 0.1  & 0.41 & 0.71 & 0.44 \\
        IRC            & 0.19 & 0.25 & 0.3  & 0.23 & 0.53 & 0.13 & 0.19  \\
        ITC            & 0.37 & 0.44 & 0.5  & 0.37 & 0.73 & 0.57 & 0.52 \\
        \bottomrule
    \end{tabular}
\caption{Item-level characteristics of the standard step-by-step test}
    \label{tb_IRC_ITC_SbS}
\end{table}

The mean scores were generally high, indicating that most solution steps were relatively easy for the analysed cohort. This result is consistent with the concentration of total scores in the upper part of the scale reported above.

The discrimination indices separate the solution steps into two groups. The initial steps (\textbf{u}, \textbf{dv}, \textbf{du}, and \textbf{v}) showed limited discrimination, with values ranging from $0.09$ to $0.26$. Among these steps, \textbf{du} approached the commonly used threshold of $0.3$, whereas \textbf{u} and \textbf{v} contributed very little to differentiating between students. In contrast, the later steps provided stronger discrimination: \textbf{uv} reached $0.41$, \textbf{vdu} reached $0.71$, and the final integration step (\textbf{int}) reached $0.44$.

The correlation-based indicators provide a more nuanced picture. All ITC values exceeded $0.3$, indicating that each response field was positively related to the total test score. However, the IRC values were more uneven. The \textbf{uv} step showed the most balanced psychometric profile, combining good discrimination ($D=0.41$), the highest IRC value ($0.53$), and the strongest ITC value ($0.73$). By contrast, \textbf{vdu} had the highest discrimination index ($D=0.71$) but a low IRC value ($0.13$), suggesting that this step separated high- and low-performing students while being weakly aligned with the remaining solution steps when considered independently. A similar pattern was observed for the final integration step, which showed acceptable discrimination ($D=0.44$) but a low IRC value ($0.19$).

The correlation matrix in Table~\ref{tb_corr_mtrx_SbS} provides additional insight into the internal structure of the test. The relationships among solution steps were uneven: some adjacent or conceptually related steps showed moderate correlations, whereas several correlations involving the early and final steps were close to zero. This pattern is expected for a step-by-step task, because the response fields represent consecutive stages of a single solution rather than independent test items. Low variability in the initial steps and the possible propagation of errors across later stages may both weaken the inter-item correlations.

\begin{table}[H]
    \centering
    \includegraphics[scale=0.85]{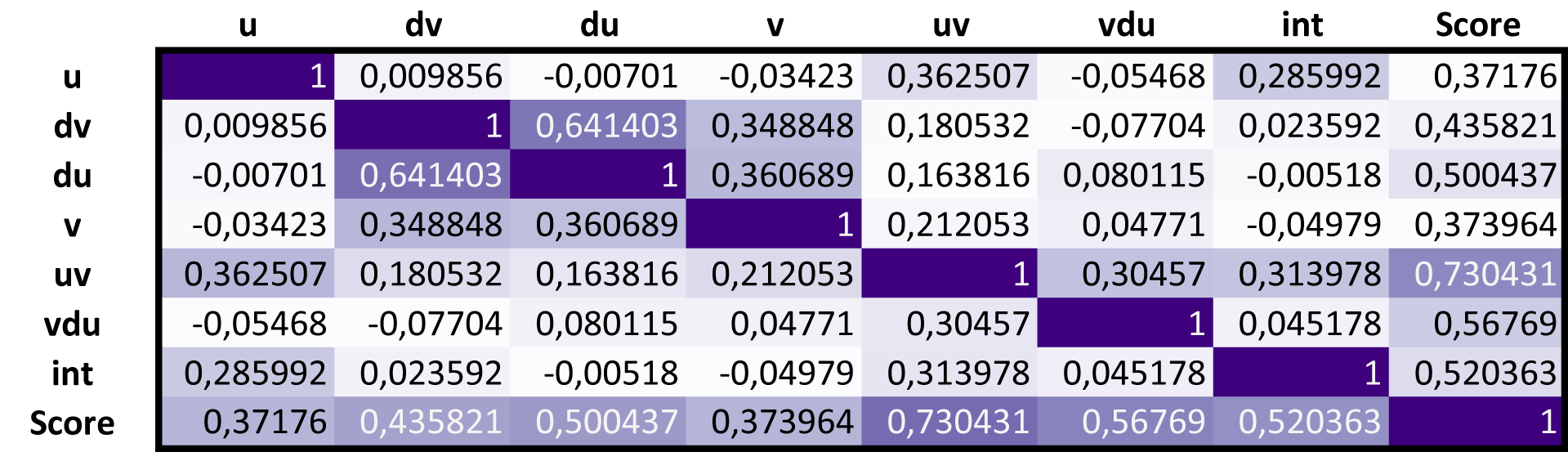}
    \caption{Heatmap of correlation matrix of item scores}
    \label{tb_corr_mtrx_SbS}
\end{table}

These item-level results show that the standard step-by-step test captures the main stages of the solution process, but its diagnostic contribution is uneven across steps. The initial fields mainly verify basic procedural choices, whereas the later fields, especially \textbf{uv}, \textbf{vdu}, and \textbf{int}, provide most of the differentiation between students. This should be considered when comparing the standard format with the STACK-based implementation in the following subsection.

\subsection{Analysis of the STACK-Based Step-by-Step Test}

The STACK-based step-by-step test used the same response-field structure as the standard test and applied the symbolic validation approach described in Section~\ref{subsec:STACK_SbST&SV}. In the analysed administration, the parameters were sampled from predefined sets, $A\in\{2,3,\ldots,9\}$, $B\in\{-5,\ldots,5\}$, $C\in\{2,3,4,5\}$
to generate integrands of the form~ \eqref{undr_int_fn}. Although the parameter space allowed $704$ possible variants, $30$ generated variants were validated and used in the present assessment. Each student received one of these variants during the test.

The implementation distinguished full- and partial-credit cases through separate rules. No post-hoc expansion of enumerated answer forms was needed, although two isolated item scores were corrected after manual review. One additional input-related case could not be evaluated automatically and was also reviewed manually. These decisions were incorporated into the final scoring matrix used in the subsequent analyses.

As shown in Table~\ref{tb_STACK_NoV}, most solution steps required only one validation rule. The steps \textbf{v} and \textbf{int} required two rules because the scoring scheme distinguished between different credit levels or acceptable mathematical forms. In the final integration step, for example, one rule corresponded to the fully correct antiderivative including the arbitrary constant, whereas another assigned partial credit to an otherwise correct antiderivative written without the constant. Thus, unlike the standard test, the number of validation rules did not increase in proportion to the variety of equivalent student representations.

\begin{table}[h]
\centering
\begin{tabular}{l|ccccccc}
\toprule
Characteristics  &  u & dv & du & v & uv & vdu & int \\
\midrule
Validation rules  & 1 & 1 & 1 & 2 & 1 & 1 & 2  \\
\bottomrule
\end{tabular}
\caption{Number of STACK validation rules specified for each solution step}
\label{tb_STACK_NoV}
\end{table}

The next part examines whether these implementation features are reflected in the distribution of total scores and the internal-consistency estimates of the STACK-based test.

\begin{figure}[H]
  \center{\includegraphics[scale=0.25]{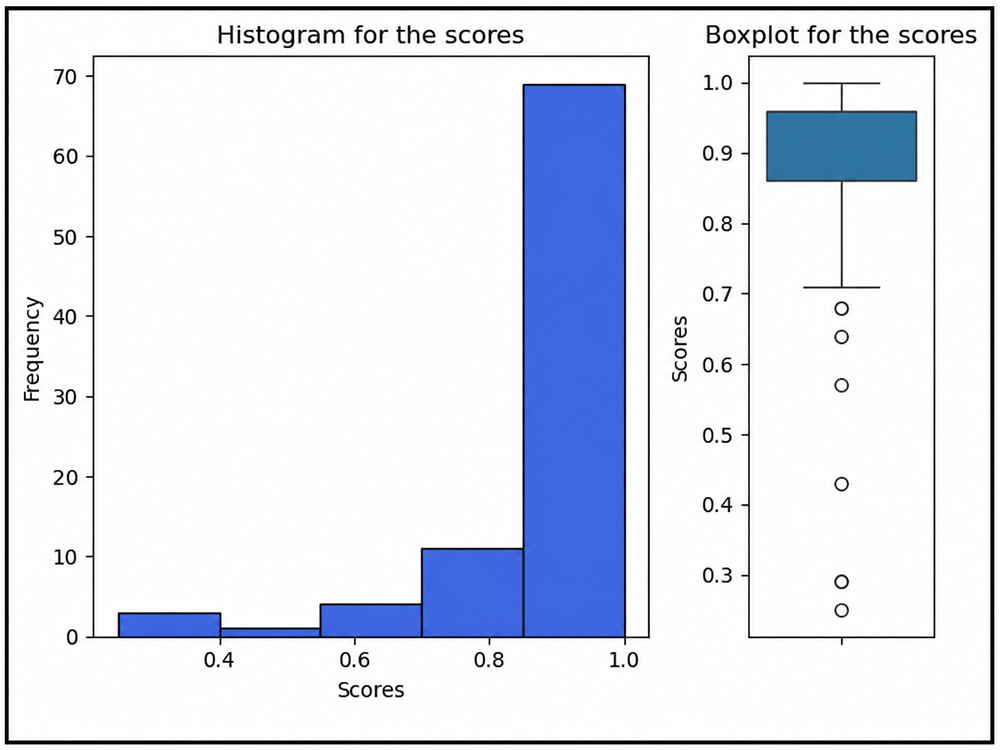}}
  \caption{Histogram and boxplot of the results for the STACK-based step-by-step test}
  \label{fig_hist_STACK}
 \end{figure}

Based on the final manually verified scoring matrix, Figure~\ref{fig_hist_STACK} summarises the distribution of total scores in the STACK-based step-by-step test. The histogram shows that most scores are concentrated in the upper part of the scale, while the boxplot reveals a small number of lower-scoring observations. These observations were retained in the analysis because they represent plausible variation in student performance rather than data errors.

Table~\ref{tb_stat_char_STACK} confirms this pattern. The mean score was $0.886$, while both the median and the mode were equal to $0.96$, indicating that most students achieved high overall results. The strong negative skewness ($-2.5$) and high kurtosis ($6.5$) indicate a distribution sharply concentrated near the upper end of the scale, with a small number of substantially lower scores.

\begin{table}[h]
    \centering
    \begin{tabular}{cccccc}
        \toprule
         \textbf{Mean} & \textbf{Median} & \textbf{Mode} &\textbf{SD} & \textbf{Skewness} & \textbf{Kurtosis} \\
        \midrule
        0.886 & 0.96 & 0.96 & 0.156 & -2.5 & 6.5 \\
        \bottomrule
    \end{tabular}

    \vspace*{1mm}
    \begin{tabular}{cc}
        \toprule
         \textbf{Cronbach's $\alpha$} & \textbf{McDonald's $\omega$}\\
        \midrule
         0.69 & 0.73 \\
        \bottomrule
    \end{tabular}
    \caption{Descriptive statistics and reliability coefficients for the STACK-based test}
    \label{tb_stat_char_STACK}
\end{table}

The concentration of scores near the maximum suggests a pronounced ceiling effect. As a result, the test captured overall procedural success well but had limited capacity to differentiate among students at the upper end of the performance distribution.

The reliability estimates provide a more favourable picture of internal consistency. Cronbach's $\alpha$ was $0.69$, close to the commonly used threshold of $0.7$, while McDonald's $\omega$ reached $0.73$, indicating acceptable internal consistency. This pattern is consistent with the more structured scoring mechanism of STACK, where mathematically equivalent expressions and partially correct responses can be handled through symbolic validation rules.

The STACK-based test therefore provided a relatively coherent overall score structure, although the high success rates continued to limit differentiation among stronger students. The item-level analysis below examines which solution steps contributed most to this pattern.

Item-level characteristics of the STACK-based test were examined for each solution step using mean scores, standard deviations, discrimination indices, item-rest correlations (IRC), and item-total correlations (ITC), as reported in Table~\ref{tb_IRC_ITC_STACK}.

\begin{table}[h]
    \centering
    \begin{tabular}{c|ccccccc}
        \toprule
        Response field  &   u  &  dv  &  du  &  v   &  uv  &  vdu & int \\
        \midrule
        Mean Score     & 0.96 & 0.95 & 0.95 & 0.95 & 0.88 & 0.68 & 0.86 \\
        SD             & 0.21 & 0.2  & 0.2 & 0.21 & 0.33  & 0.35 & 0.31 \\
        Discrimination & 0.17 & 0.17 & 0.19 & 0.17  & 0.44 & 0.52 & 0.41 \\
        IRC            & 0.26 & 0.4 & 0.51  & 0.68 & 0.56 & 0.33 & 0.21  \\
        ITC            & 0.43 & 0.55 & 0.64  & 0.78 & 0.75 & 0.6 & 0.46 \\
        \bottomrule
    \end{tabular}
    \caption{Item-level characteristics of the STACK-based step-by-step test}
    \label{tb_IRC_ITC_STACK}
\end{table}

Mean scores ranged from $0.68$ to $0.96$, showing that the initial steps were completed successfully by most students, whereas later steps introduced greater difficulty. The first four steps (\textbf{u}, \textbf{dv}, \textbf{du} and \textbf{v}) had very high mean scores ($0.95$-$0.96$) and low discrimination indices ($0.17$-$0.19$). This indicates that these steps verified basic procedural components of the solution but contributed relatively little to differentiating between students.

Although the first four steps had similarly low discrimination indices, their correlation-based indicators differed. In particular, \textbf{du} and especially \textbf{v} showed stronger associations with the rest of the test. The step \textbf{v} had the highest IRC ($0.68$) and ITC ($0.78$) values, suggesting that successful completion of this stage was closely aligned with overall test performance, even though it was not difficult enough to provide strong discrimination.

The later solution steps provided the main source of differentiation between students. The \textbf{uv} step showed the most balanced psychometric profile, combining good discrimination ($D=0.44$), a strong IRC value ($0.56$), and a high ITC value ($0.75$). The \textbf{vdu} step was the most difficult and the most discriminative, with the lowest mean score ($0.68$) and the highest discrimination index ($D=0.52$). The final integration step also contributed to differentiation ($D=0.41$), although its lower IRC value ($0.21$) suggests that performance at this stage may depend on accumulated errors and on the partial-credit structure of the final answer.

The item-level results therefore show that the STACK-based test combines high success rates on the initial procedural steps with stronger diagnostic contribution from the later stages of the solution. The generally positive IRC and ITC values are consistent with a coherent relationship between individual solution steps and overall test performance.

The correlation matrix in Table~9 provides additional insight into the internal structure of the assessment. The intermediate and later steps show a more coherent pattern of associations than the initial steps. In particular, the steps related to the construction and application of the integration-by-parts formula (\textbf{v}, \textbf{uv}, and \textbf{vdu}) are more consistently associated with the remaining solution process than the initial selection of \textbf{u} and \textbf{dv}. At the same time, several weak or near-zero correlations remain, especially for steps with very high success rates or for the final response field. This pattern is expected for a step-by-step task, because the response fields represent consecutive stages of a single solution rather than independent test items. Low variability in the initial steps and the possible accumulation of errors across later stages may both affect the inter-item correlations.

\begin{table}[H]
    \centering
    \includegraphics[scale=0.85]{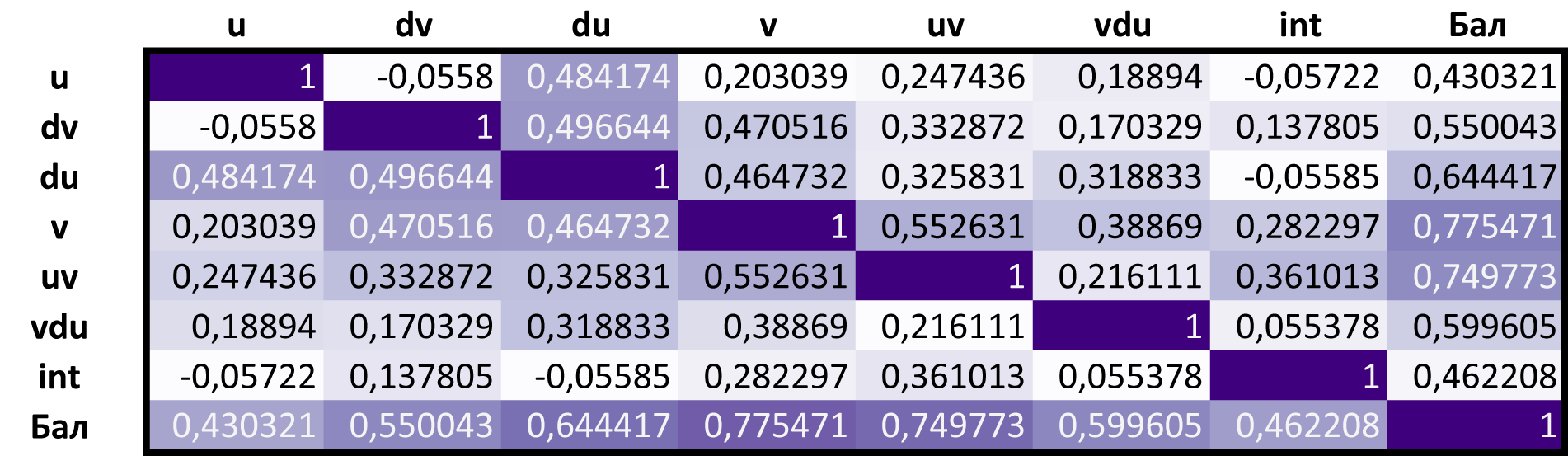}
    \caption{Heatmap of correlation matrix of item scores}
    \label{tb_corr_mtrx_STACK}
\end{table}

These results indicate that the STACK-based test captures the logical structure of the integration-by-parts solution while assigning most diagnostic value to the later, more demanding steps. They also provide the basis for the direct comparison with the standard step-by-step test in the following subsection.

\subsection{Comparative Analysis of the Standard and STACK-Based Step-by-Step Tests}

The preceding subsections analysed the two tests separately. This subsection compares the empirical evidence for the implementation differences described in Section~\ref{subsec:Comp_of_VM&IC} and then examines whether these differences are reflected in the psychometric indicators.

The standard test contained $16$ manually prepared variants, whereas the STACK-based implementation used parameterised templates with a parameter space of $704$ possible structurally comparable variants, $30$ of which were validated and used in the present assessment. These figures characterise the difference in variant-generation workload between the two implementations.

The implementation data also confirm the contrast in scoring mechanisms described in Section~\ref{subsec:Comp_of_VM&IC}. In the standard test, the set of predefined scoring patterns had been cumulatively expanded by $29$ entries relative to the initial version. This expansion reflected alternative full-credit representations and common partially correct forms of student input. In the STACK-based test, the same role was handled by a small set of symbolic validation rules, without post-hoc expansion of enumerated answer forms.

The scoring data also show a clear difference in post-test moderation. In the present administration, the standard test required\ \, $25$\ \, manual score corrections, whereas the STACK-based test required only two. In addition, four responses in the standard test and one response in the STACK-based test could not be evaluated automatically because of input-related issues. All these cases were reviewed manually and incorporated into the final scoring matrices used for the psychometric analyses.

Taken together, these implementation indicators show that the two formats differed not only in their validation mechanisms but also in the amount of post-test moderation required. The next part examines whether these differences were accompanied by differences in score distributions, reliability, and item-level indicators.

Both tests produced high overall scores, confirming that most students successfully completed the integration-by-parts task. The mean score was $0.851$ for the standard test and $0.886$ for the STACK-based test, indicating broadly similar levels of overall performance (Table~\ref{tb_stat_char_both}).

The distributions differed more noticeably in shape. The STACK-based assessment showed stronger negative skewness ($-2.50$ compared with $-1.25$) and substantially higher kurtosis ($6.50$ compared with $1.28$), indicating a stronger concentration of scores near the upper end of the scale. The interquartile interval for the STACK-based test was also shifted upward ($Q_1=0.88$, $Q_3=0.96$) relative to the standard test ($Q_1=0.821$, $Q_3=0.964$). At the same time, the STACK-based test had a narrower interquartile range, reflecting a stronger clustering of scores among high-performing students. These results indicate that a ceiling tendency was present in both formats but was more pronounced in the STACK-based assessment.

\begin{table}[h]
\centering
\begin{tabular}{lcccccc}
\toprule
\textbf{Test type} & \textbf{Mean} & \textbf{SD} & \textbf{Skewness} & \textbf{Kurtosis} & \textbf{$\alpha$} & \textbf{$\omega$} \\
\midrule
Standard & 0.851 & 0.142 & -1.25 & 1.28 & 0.49 & 0.59 \\
STACK-Based & 0.886 & 0.156 & -2.5 & 6.5  & 0.69 & 0.73 \\
\bottomrule
\end{tabular}
\caption{Comparison of descriptive statistics and reliability coefficients}
\label{tb_stat_char_both}
\end{table}

The reliability coefficients show a clearer difference between the two formats. For the standard test, Cronbach's $\alpha$ was low ($0.49$), and McDonald's $\omega$ remained below the commonly used threshold for acceptable internal consistency ($0.59$). In the STACK-based test, Cronbach's $\alpha$ increased to $0.69$, approaching the conventional threshold of $0.7$, while McDonald's $\omega$ reached $0.73$, indicating acceptable internal consistency.

The higher observed reliability coefficients are consistent with the STACK scoring mechanism, in which equivalent full-credit responses and partially correct responses are handled through explicit symbolic validation rules rather than through cumulatively expanded predefined scoring patterns. Thus, the STACK-based format demonstrated stronger internal consistency, although the high concentration of scores near the maximum continued to limit differentiation among higher-performing students.

Table~\ref{tb_IRC_ITC_both} shows that both assessment formats exhibit a similar overall profile across solution steps: the initial steps were completed successfully by most students, whereas the later steps provided most of the differentiation between performance levels.

\begin{table}[h]
\centering
\begin{tabular}{cc|ccccccc}
\toprule
\textbf{Characteristic} & \textbf{Test}&  u & dv & du & v & uv & vdu & int \\
\midrule
\multirow{2}{*}{Mean Score} & Standard & 0.96 & 0.93 & 0.91 & 0.97 & 0.89 & 0.54 & 0.76 \\
 & STACK &  0.96 & 0.95 & 0.95 & 0.95 & 0.88 & 0.68 & 0.86 \\
\midrule
\multirow{2}{*}{Discrimination} & Standard &  0.09 & 0.21 & 0.26 & 0.10 & 0.41 & 0.71 & 0.44 \\
 & STACK &  0.17 & 0.17 & 0.19 & 0.17 & 0.44 & 0.52 & 0.41 \\
\bottomrule
\multirow{2}{*}{IRC} & Standard &  0.19 & 0.25 & 0.3  & 0.23 & 0.53 & 0.13 & 0.19  \\
  & STACK &  0.26 & 0.4 & 0.51  & 0.68 & 0.56 & 0.33 & 0.21  \\
\bottomrule
\multirow{2}{*}{ITC} & Standard &  0.37 & 0.44 & 0.5  & 0.37 & 0.73 & 0.57 & 0.52 \\
 & STACK & 0.43 & 0.55 & 0.64  & 0.78 & 0.75 & 0.6 & 0.46 \\
\bottomrule
\end{tabular}
\caption{Comparison of item-level characteristics for both tests}
 \label{tb_IRC_ITC_both}
\end{table}

Mean scores were broadly comparable across the two formats. The most noticeable differences occurred in the later stages, where the STACK-based test produced higher mean scores for \textbf{vdu} and \textbf{int}. This pattern is consistent with the reduced dependence of STACK scoring on predefined input patterns and with its ability to assign full or partial credit through symbolic validation rules.

The discrimination indices show the same general tendency in both formats. The initial steps had limited discriminatory power, whereas the later steps contributed more strongly to differentiating between students. The standard test showed the highest discrimination for \textbf{vdu} ($D=0.71$), while the STACK-based test showed a more even distribution of discrimination across the later steps ($D=0.44$ for \textbf{uv}, $D=0.52$ for \textbf{vdu}, and $D=0.41$ for \textbf{int}).

The strongest difference between the two formats is observed in the correlation-based indicators. The STACK-based test produced higher IRC values for all solution steps and higher ITC values for most of them. This indicates a stronger alignment between individual response fields and the overall score structure.

The difference is particularly clear for the step \textbf{v}. In the STACK-based test, this step had the highest IRC ($0.68$) and ITC ($0.78$), whereas the corresponding values in the standard test were substantially lower ($0.23$ and $0.37$, respectively). Since this step represents the transition from selecting the integration components to applying the integration-by-parts formula, the result suggests that the STACK-based scoring scheme captured this intermediate stage more consistently.

The final integration step (\textbf{int}) requires a more cautious interpretation. Although the STACK-based test produced a higher mean score for this step, its ITC value was slightly lower than in the standard test. This may reflect the partial-credit structure of the final answer, where an otherwise correct antiderivative without the arbitrary constant was distinguished from a fully correct response.

The comparative item-level evidence therefore suggests that the advantage of the STACK-based test lies not in uniformly higher discrimination across all steps, but in a more coherent relationship between individual solution steps and the total score. This supports the interpretation that symbolic validation provides a more stable scoring structure while preserving the step-by-step nature of the assessment.

The practical implications of these differences for test development and maintenance are summarised in the next subsection.

\subsection{Practical Recommendations}

The results of the comparative analysis indicate that the choice between the two assessment formats should depend on the purpose, scale, expected reuse of the test, and the required level of scoring robustness.

For small or infrequently used assessments, conventional step-by-step tests remain a feasible solution. For such use cases, their simpler construction may outweigh the need for occasional manual review. However, their use becomes less efficient when many task variants are needed or when student responses include diverse full- and partial-credit forms. In such cases, the cumulative maintenance of predefined scoring patterns may require substantial manual review.

For large cohorts, repeated administrations, or assessments involving several parameterised variants, the STACK-based format is generally more appropriate. Its symbolic validation rules reduce dependence on enumerated input patterns, support full- and partial-credit scoring, and allow comparable variants to be generated from parameterised templates. Although the initial development of STACK tasks requires greater technical and mathematical effort, this cost is offset when tests are reused, administered to larger cohorts, or generated from parameterised templates.

The findings also suggest that STACK is particularly suitable for multi-step mathematical problems where the structure of the solution is part of what should be assessed. In such tasks, symbolic validation can support a closer correspondence between automated scoring and the mathematical meaning of students' responses.

Finally, the ceiling effect observed in both formats should be addressed in future test development. Later solution steps or additional extension tasks should be made more diagnostically demanding in order to improve differentiation among higher-performing students.

\section{Discussion and Conclusions}
\label{sec:Disc&Conc}

The comparative analyses presented in Section~\ref{sec:CompAn42Tests} show that the two assessment formats differed primarily in their scoring mechanisms rather than in the mathematical structure of the task. Both tests assessed the same integration-by-parts procedure and produced high overall student performance. The main distinction concerned how full-credit and partial-credit responses were represented and evaluated: the standard test relied on cumulatively expanded predefined scoring patterns, whereas the STACK-based test used symbolic validation rules.

The STACK-based assessment showed a more coherent score structure, as reflected in higher reliability estimates and stronger item-rest relationships across solution steps. The difference was especially visible for the intermediate step corresponding to the computation of \textbf{v}, which showed the strongest alignment with the overall score structure in the STACK-based test. This suggests that symbolic validation captured this stage of the integration-by-parts procedure more consistently than predefined input matching.

These results are consistent with earlier work on computer-aided mathematics assessment, which emphasises that symbolic systems improve robustness by evaluating mathematical meaning rather than surface textual form~\cite{sangwin2013,SangwinKocher2016}. They also support broader arguments that effective assessment of mathematics should capture reasoning processes and not only final answers~\cite{beevers2003,greenhow2015}. The contribution of the present study is to show that this distinction is reflected not only in grading practice but also in psychometric indicators, particularly reliability estimates and item-rest relationships.

The practical implication is a trade-off between initial development effort and long-term scoring stability. Standard step-by-step tests are easier to implement initially, but their dependence on cumulatively expanded scoring patterns increases the need for manual moderation. STACK requires greater initial expertise, especially in designing symbolic validation rules and answer trees, but it provides a more scalable structure for repeated and parameterised assessments.

A limitation common to both formats was the concentration of scores near the upper end of the scale. This pattern suggests that the selected integration-by-parts task verified procedural competence effectively but was not sufficiently demanding to differentiate among higher-performing students. The ceiling effect therefore appears to be associated primarily with the characteristics of the assessment task rather than with symbolic validation itself.

Several limitations should be considered when interpreting the findings. First, the comparison involved two existing student cohorts rather than randomly assigned groups, so differences in prior mathematical preparation cannot be ruled out. Second, the study focused on a single mathematical topic, integration by parts, within one undergraduate course at a single institution. Third, the psychometric analyses were based on final manually verified scoring matrices. This ensured fair scoring of mathematically meaningful responses but means that the reported reliability estimates describe the assessments after manual review rather than purely automated scoring. Finally, although Classical Test Theory provided an appropriate framework for this comparison, complementary psychometric approaches may provide additional insight into the behaviour of symbolic assessment items.

Future research should examine whether the observed psychometric differences remain stable across other mathematical topics involving different forms of symbolic reasoning. Extending the comparison to algebraic manipulation, differential equations, and multivariable calculus would help determine the generalisability of the present findings. Studies involving larger and more heterogeneous student populations would also enable analyses based on Item Response Theory and Differential Item Functioning. Rasch-based approaches could additionally be used to examine item functioning, measurement invariance, and differences between cohorts or task variants~\cite{zeileis2025}. Another useful direction would be to compare automated raw scores with manually verified scores in order to quantify the effect of post-test moderation on reliability and item-level indicators. Finally, symbolic assessment could be extended beyond routine procedural exercises toward more sophisticated mathematical reasoning tasks~\cite{rowlett2022,SangwinBickerton2021}.

The study shows that the main value of STACK-based assessment lies not merely in automating grading but in providing a scoring approach that is more closely aligned with the symbolic nature of mathematical reasoning. By evaluating mathematical equivalence and supporting structured full- and partial-credit scoring, symbolic validation offers a stronger basis for reliable, scalable, and pedagogically meaningful assessment of multi-step problems in university-level mathematics and related STEM disciplines.

\section*{Data Availability}
The data used in this study were collected following protocols that ensure the anonymity of participants. Names and other identifying information have not been published and cannot be shared externally. However, interested researchers may contact the authors for additional details regarding the analysis methods or aggregated results.

\section*{Conflicts of Interest}
The authors declare no conflicts of interest related to this study. All work was conducted
independently and without commercial influence.


\end{document}